\documentclass[12pt]{article}
\usepackage[latin1]{inputenc}
\usepackage[T1]{fontenc}
\usepackage{amsmath,amssymb}
\usepackage[all]{xy}
\usepackage{amsfonts}
\usepackage{amssymb}
\usepackage{lambda}
\usepackage[oldstyle]{kpfonts}
\def\tab{\hspace*{6mm}}

\newcounter{numero}
\newcommand{\num}
{\refstepcounter{numero}
\noindent \thenumero.---}

\newcommand{\fin}{
\vskip 2mm
\noindent
$\Box$}

\begin{document}
\begin{center}
\LARGE Le Probl\`eme Inverse de Galois sur les corps des fractions tordus \`a ind\'etermin\'ee centrale
\end{center}
\vskip 3mm
\noindent
\centerline{Bruno DESCHAMPS et Fran\c cois LEGRAND}
\vskip 6mm
\noindent
\begin{center}
\begin{minipage}{11cm}
{\small {\bf Abstract.---} In this article, we show that the Inverse Galois Problem over a skew field $H$ of finite dimension over its center $k$ is equivalent to a variant of the Inverse Galois Problem over $k$ involving a polynomial constraint. As an application, we show that if $k$ contains an ample field, then the Inverse Galois Problem has a positive answer over the skew field $H(t)$ of rational fractions with central indeterminate.\\
{\bf R\'esum\'e.---} Dans cet article, nous montrons que le Probl\`eme Inverse de Galois sur un corps gauche $H$ de dimension finie sur son centre $k$ est \'equivalent \`a une variante du Probl\`eme Inverse de Galois sur $k$ faisant intervenir une contrainte polynomiale. En application de ce r\'esultat nous montrons que, si $k$ contient un corps ample, alors le Probl\`eme Inverse de Galois admet une r\'eponse positive sur le corps $H(t)$ des fractions rationnelles tordu \`a ind\'etermin\'ee centrale.}
\end{minipage}
\end{center}
\vskip 1cm
\noindent
\tab Le probl\`eme inverse de la th\'eorie de Galois sur un corps $K$ s'\'enonce de la mani\`ere suivante :
\vskip 2mm
\noindent
{\bf $\hbox{\bf PIG}_K$ (Probl\`eme Inverse de Galois sur $K$)} : {\it Tout groupe fini $G$ est-il le groupe de Galois d'une extension galoisienne $L/K$ ?}
\vskip 2mm
\noindent
Etant donn\'ee une extension de corps commutatifs $H/k$, il est clair que, de mani\`ere g\'en\'erale, l'implication $\hbox{\rm PIG}_H\Longrightarrow \hbox{\rm PIG}_k$ est fausse, m\^eme lorsque $H/k$ est finie. On peut, par exemple, penser au cas du corps $k={\mathbb Q}^{\hbox{\rm \scriptsize tr}}$ des nombres alg\'ebriques totalement r\'eels, corps sur lequel seuls les groupes finis engendr\'es par des involutions sont groupes et Galois (voir \cite{FHV93}) et $H={\mathbb Q}^{\hbox{\rm \scriptsize tr}}(\sqrt{-1})$ dont le groupe de Galois absolu est prolibre (voir, par exemple, \cite[Remark 7.10(b)]{GJ02}). L'\'etude de la r\'eciproque $\hbox{\rm PIG}_k\Longrightarrow \hbox{\rm PIG}_H$ est plus commode : si l'on se donne une extension galoisienne $L/k$ de groupe fini $G$ telle que l'alg\`ebre $M=L\otimes_k H$ soit un corps, alors un argument classique de th\'eorie de Galois montre que $M/H$ est une extension galoisienne de groupe $G$. Comme $L/k$ est une extension alg\'ebrique, la condition "$L\otimes_k H$ est un corps" \'equivaut juste au fait que les extensions $L/k$ et $H/k$ sont lin\'eairement disjointes sur $k$. Si l'on suppose l'extension $H/k$ finie alors la lin\'eaire disjonction peut s'interpr\'eter de la mani\`ere suivante : on consid\`ere $(e_1,\dots ,e_n)$, une $k$-base de $H$ et l'on consid\`ere la forme
$$P_H (x_1,\dots ,x_n)=N_{H/k}(x_1e_1+\cdots +x_ne_n)\in k[x_1,\dots ,x_n]$$
o\`u $N_{H/k}$ d\'esigne la norme de l'extension $H/k$. Les extensions $H/k$ et $L/k$ sont alors lin\'eairement disjointes sur $k$ si et seulement si la forme $P_H$ ne poss\`ede que le z\'ero trivial sur $L$. On peut donc consid\'erer une variante du Probl\`eme Inverse de Galois :
\vskip 2mm
\noindent
{\bf $\hbox{\bf PIG${\cal F}$C}_k$ (Probl\`eme Inverse de Galois \`a ${\cal F}$-contrainte)} : {\it Pour la $k$-forme ${\cal F}$ donn\'ee, tout groupe fini $G$ est-il le groupe de Galois d'une extension galoisienne $L/k$ telle que ${\cal F}$ ne poss\`ede que le z\'ero trivial sur $L$ ?}
\vskip 2mm
\noindent
Pour le choix de ${\cal F}=P_H$, on voit donc que $\hbox{\rm PIG${\cal F}$C}_k\Longrightarrow \hbox{\rm PIG}_H$. En fait, un argument classique montre que, pour ce choix de ${\cal F}$, on a $\hbox{\rm PIG${\cal F}$C}_k\Longleftrightarrow \hbox{\rm PIG}_k$\footnote{Notons $n$ le nombre de corps interm\'ediaires de $H/k$. Pour un groupes $G$ donn\'e on consid\`ere une extension galoisienne $L/k$ de groupe $G^{n+1}$ et l'on consid\`ere le corps $L_i$ des invariants de $L$ par le sous-groupe de  $G^{n+1}$ obtenu en retirant le i-\`eme facteur direct. Les extensions $L_i/k$ sont alors galoisiennes de groupe $G$ et elles sont lin\'eairement disjointes deux \`a deux sur $k$. Puisqu'elles sont au nombre de $n+1$, au moins l'une d'entre elle est d'intersection avec $H$ r\'eduite \`a $k$, c'est-\`a-dire lin\'eairement disjointe de $H$ sur $k$.}, mais c'est la formulation $\hbox{\rm PIG${\cal F}$C}_k$ qui va nous int\'eresser dans cet article.
\vskip 2mm
\noindent
\tab Traditionnellement, les arithm\'eticiens consid\`erent le(s) probl\`eme(s) inverse(s) de Galois pour des corps commutatifs. Nonobstant, la th\'eorie de Galois poss\`ede une g\'en\'eralisation au cas non commutatif et il est donc possible de regarder le Probl\`eme Inverse de Galois dans le cas des corps gauches. La d\'efinition la plus g\'en\'erale d'extension galoisienne est celle donn\'ee par Artin : une extension $L/k$ est dite galoisienne si le corps des invariants de $L$ sous l'action de $\hbox{\rm Aut}(L/k)$ est \'egal \`a $k$. C'est avec cette d\'efinition que s'entend la th\'eorie de Galois des corps gauches. Le r\'esultat central de cet article est le
\vskip 2mm
\noindent
{\bf Th\'eor\`eme A.---} {\it Soient $k$ un corps commutatif et $H$ un corps de dimension finie sur son centre $k$. Si ${\cal F}$ d\'esigne la $k$-forme associ\'ee \`a la norme r\'eduite de l'extension $H/k$ relativement au choix d'une $k$-base de $H$, alors on a l'\'equivalence
$$\hbox{\rm PIG${\cal F}$C}_k\Longleftrightarrow \hbox{\rm PIG}_{H}$$}
\vskip 2mm
\noindent
\tab Ainsi, en jouant sur les analogues "norme/norme r\'eduite" pour la forme ${\cal F}$ et "commutatif/gauche de centre $k$" pour le corps $H$, on voit que l'on a toujours
$$\hbox{\rm PIG${\cal F}$C}_k\Longrightarrow \hbox{\rm PIG}_{H}$$
mais, chose surprenante, l'implication $\hbox{\rm PIG}_H\Longrightarrow \hbox{\rm PIG${\cal F}$C}_k$ n'est vraie en toute g\'en\'eralit\'e que dans le cas non commutatif, comme nous l'avons illustr\'e tout \`a l'heure. 
\vskip 2mm
\noindent
\tab La premi\`ere partie de cet article est consacr\'ee \`a la preuve du th\'eor\`eme A. Nous commen\c cons cette partie par un bref rappel des \'el\'ements fondamentaux de la g\'en\'eralisation de la th\'eorie de Galois au cas des corps gauches afin de faciliter la compr\'ehension de ce texte aux arithm\'eticiens non familiaris\'es \`a cette situation. 
\vskip 2mm
\noindent
\tab La deuxi\`eme partie de cet article s'int\'eresse plus particuli\`erement au cas de certains corps de fractions rationnelles tordus. Rappelons que, si $H$ d\'esigne un corps (commutatif ou non), on d\'efinit {\it l'anneau  des polyn\^omes tordu $H[t]$ \`a ind\'etermin\'ee centrale et \`a coefficients dans $H$} comme le $H$-espace vectoriel de base $\{t^n\}_{n\geq 0}$ sur lequel on consid\`ere le produit qui v\'erifie $at=ta$ pour tout $a\in H$. Il s'agit en fait, pour la construction de Ore (voir \cite{Ore33}), de l'anneau de polyn\^omes tordu $K[t,\alpha,\delta]$ obtenu en prenant pour $\delta$ la d\'erivation nulle et pour $\alpha$ l'endomorphisme identit\'e. Cet anneau est un anneau de Ore et, en particulier, il poss\`ede un unique corps de fractions que l'on note $H(t)$ et que l'on appelle {\it le corps des fractions rationnelles tordu \`a ind\'etermin\'ee centrale et \`a coefficients dans $H$}. En tant que corps des fractions d'un anneau de Ore, le corps $H(t)$ jouit d'une propri\'et\'e tr\`es forte : tous ses \'el\'ements peuvent s'\'ecrire sous la forme $p(t)q(t)^{-1}$ avec $p(t), \, q(t)\in H[t]$. Lorsque $H$ est commutatif, le corps $H(t)$ obtenu est alors le classique corps des fractions rationnelles \`a coefficients dans $H$. Le th\'eor\`eme A montre alors que $\hbox{\rm PIG${\cal F}$C}_{k(t)}\Longleftrightarrow \hbox{\rm PIG}_{H(t)}$ pour le choix de la forme ${\cal F}$ associ\'ee \`a la norme r\'eduite de l'extension $H/k$. 
\vskip 2mm
\noindent
\tab Il existe une famille importante de corps commutatifs $k$ pour lesquels le $\hbox{\rm PIG${\cal F}$C}_{k(t)}$ admet une r\'eponse positive pour toute $k$-forme ${\cal F}$ ne poss\`edant pas de z\'ero sur le corps $k(t)$ : c'est celle des corps amples. Il s'agit d'une propri\'et\'e introduite par Pop dans \cite{Pop96} et qui joue un r\^ole majeur en th\'eorie inverse de Galois. Un corps commutatif $k$ est dit {\it ample} si toute courbe lisse d\'efinie sur $k$ et g\'eom\'etriquement irr\'eductible poss\`ede une infinit\'e de points $k$-rationnels d\`es qu'elle en poss\`ede au moins un. Parmi ces corps, on trouve les corps PAC\footnote{Un corps $k$ est dit Pseudo Alg\'ebriquement Clos (PAC) si toute $k$-vari\'et\'e non vide et g\'eom\'etriquement irr\'eductible a un ensemble Zariski dense de points $k$-rationnels.} (e.g. les corps s\'eparablement clos, le corps ${\mathbb Q}^{\hbox{\rm \scriptsize tr}}(\sqrt{-1})$), les corps valu\'es complets (e.g. ${\mathbb R}$, ${\mathbb Q}_p$, $\kappa((x))$), le corps ${\mathbb Q}^{\hbox{\rm \scriptsize tp}}$ des nombres alg\'ebriques totalement $p$-adiques, le corps ${\mathbb Q}^{\hbox{\rm \scriptsize tr}}$ des nombres alg\'ebriques totalement r\'eels (voir les articles de survol \cite{DD97b}, \cite{BSF13} et \cite{Pop14} pour plus de d\'etails). Dans cette deuxi\`eme partie nous d\'etaillons certains liens qui existent entre cette propri\'et\'e d'amplitude et notre probl\'ematique. Ceci nous permet finalement de montrer le
\vskip 2mm
\noindent
{\bf Th\'eor\`eme B.---} {\it Si $k$ d\'esigne un corps commutatif contenant un corps ample, alors le $\hbox{\rm PIG}_{H(t)}$ admet une r\'eponse positive pour tout corps $H$ de dimension finie sur son centre $k$} 
\vskip 2mm
\noindent
\tab Une premi\`ere application de ce r\'esultat, sans doute la plus simple, s'obtient en consid\'erant $k={\mathbb R}$ : {\it tout groupe fini est groupe de Galois sur ${\mathbb H}(t)$}, o\`u ${\mathbb H}$ d\'esigne le corps des quaternions d'Hamilton. Mais le th\'eor\`eme B s'applique aussi \`a des corps ayant des gros groupes de Brauer (e.g. ${\mathbb Q}^{\hbox{\rm \scriptsize tr}}$, $\overline{k_0}(x_1,\dots ,x_n)$), ce qui montre l'\'etendue du champ d'application. 
\vskip 1cm
\noindent
{\bf \large 1.--- Preuve du th\'eor\`eme A.}
\vskip 5mm
\noindent
{\bf 1.1.--- Descente galoisienne.}
\vskip 2mm
\noindent
\tab Nous rappelons ici les r\'esultats fondamentaux de la th\'eorie de Galois. Dans le cas non commutatif, la pr\'esence d'automorphismes int\'erieurs perturbe un peu les choses. Cette perturbation se mesure de la mani\`ere suivante : si $L$ est un corps de centre $C$ et si $L/K$ d\'esigne une extension galoisienne de groupe $G$ telle que l'une des dimensions de $L$ en tant que $K$-espace vectoriel droite ou gauche soit finie, alors 
\vskip 2mm
\noindent
1/ Les dimensions droite et gauche de $L$ en tant que $K$-espace vectoriel sont \'egales (et l'on peut donc parler du degr\'e $[L:K]$ de l'extension sans se soucier de quel cot\'e l'on consid\`ere le $K$-espace vectoriel $L$).
\vskip 2mm
\noindent
2/ L'ensemble $A=\{a\in L^*\ |\ I(a)\in G\}\cup \{0\}$ est un corps (ici $I(a)$ d\'esigne l'automorphisme de conjugaison int\'erieur associ\'e \`a l'\'el\'ement $a\in L^*$). L'ensemble $A$ est en fait le commutant du corps $K$ dans $L$.
\vskip 2mm
\noindent
3/ Si l'on consid\`ere $G_0=\{I(a)\ |\ a\in A\}$ le sous-groupe de $G$ compos\'e des automorphismes int\'erieurs, alors on a 
$$[L:K]=[G:G_0][A:C]$$
\vskip 1mm
\noindent
\tab On voit qu'il peut donc exister des situations o\`u $[L:K]$ est fini et o\`u $G$ est infini (c'est par exemple le cas lorsque que l'on consid\`ere un corps $L$ de dimension finie sur son centre $K$), on peut m\^eme trouver des exemples o\`u $G$ est fini mais o\`u $|G|>[L:K]$ (voir \cite{Des18}). Pour autant, une cons\'equence de la propri\'et\'e 3/ est que l'on a $[L:K]=|G|$, d\`es que $G_0=1$ (ou de mani\`ere \'equivalente d\`es que le commutant $A$ de $K$ dans $L$ est \'egal au centre $C=Z(L)$). Dans cette situation, on dit que l'extension $L/K$ est {\it ext\'erieure}. Dans la situation oppos\'ee o\`u $G=G_0$, on dit que $L/K$ est {\it int\'erieure}.
\vskip 2mm
\noindent
\tab Venons-en maintenant aux correspondances galoisiennes. Pour les \'etudier nous allons avoir besoin de la notion de $N$-groupe. Un groupe $G$ d'automorphismes d'un corps $K$ est appel\'e $N$-groupe si l'ensemble
$$A(G,K)=\left\{x\in K^*/\ I(x)\in G\right\}\cup \{0\}$$
est un corps. Un sous-groupe $S$ d'un $N$-groupe $G$ sera dit $N$-invariant dans $G$, si le sous-$N$-groupe de $G$ engendr\'e par le normalisateur $N_G(S)$ de $S$ dans $G$ est \'egal \`a $G$ tout entier. Si $L/K$ d\'esigne une extension galoisienne de degr\'e fini alors 
\vskip 2mm
\noindent
a) Les correspondances galoisiennes op\`erent une bijection r\'eciproque l'une de l'autre entre les sous-$N$-groupes de $G=\hbox{\rm Gal}(L/K)$ et les extensions interm\'ediaires de $L/K$.
\vskip 2mm
\noindent
b) Si, par les correspondances galoisiennes, le sous-$N$-groupe $S$ correspond \`a l'extension interm\'ediaire $L_0$, alors $L/L_0$ est galoisienne de groupe $S$. De plus, le groupe $\hbox{\rm Aut}(L_0/K)$ s'identifie au groupe quotient $N_G(S)/S$ et l'extension $L_0/K$ est galoisienne si et seulement si $S$ est $N$-invariant dans $G=\hbox{\rm Gal}(L/K)$.
\vskip 2mm
\noindent
On pourra trouver dans \cite[\S3.3]{Coh95} un expos\'e complet de la th\'eorie de Galois dans le cas non commutatif et le d\'etail des preuves des propri\'et\'es que nous venons d'\'enoncer. On remarquera pour finir ce paragraphe de rappels, que dans le cas o\`u les corps sont commutatifs, on retrouve exactement les th\'eor\`emes fondamentaux de la th\'eorie de Galois usuelle.
\vskip 2mm
\noindent
\tab Si l'on consid\`ere une extension galoisienne de degr\'e fini $M/K$ de groupe de Galois $\Gamma$ et que l'on note $\Gamma_0$ l'ensemble des \'el\'ements de $\Gamma$ qui sont des automorphismes int\'erieurs, alors $\Gamma_0$ est un sous-$N$-groupe de $\Gamma$ qui est $N$-invariant (en fait normal). En posant $L=M^{\Gamma_0}$, on obtient alors la tour
$$
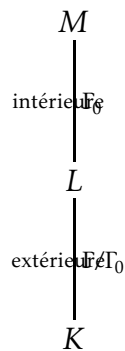
\xymatrix@!0 @R=5pc @C=5pc{\relax M\\
L \ar@{-}[u]^-*[@]{\hbox to 3pt{\hss\txt{\scriptsize int\'erieure}\hss}}_{\Gamma_0}\\
K \ar@{-}[u]^-*[@]{\hbox to 3pt{\hss\txt{\scriptsize ext\'erieure}\hss}}_{\Gamma/\Gamma_0}}$$
Ainsi, on peut toujours "d\'ecouper" une extension galoisienne de degr\'e fini $M/K$ en une partie $M/L$ int\'erieure et une partie $L/K$ ext\'erieure. Le th\'eor\`eme suivant s'int\'eresse \`a la situation duale :
\vskip 2mm
\noindent
{\bf Th\'eor\`eme\num\label{1}} {\it On suppose avoir une tour d'extensions de corps $M/L/K$ telle que $M/L$ soit galoisienne ext\'erieure de groupe fini $G$ et $L/K$ soit galoisienne int\'erieure de degr\'e fini. Si $\widetilde{\widetilde{K}}$ d\'esigne le bicommutant de $K$ dans $M$, alors l'extension $\widetilde{\widetilde{K}}/K$ est une extension galoisienne ext\'erieure de groupe $G$.}
\vskip 2mm
\noindent
{\bf Preuve :} On pose $G_0=\hbox{\rm Gal}(L/K)$, $\Gamma=\hbox{\rm Aut}(M/K)$ et l'on consid\`ere $\Gamma_0\subset \Gamma$ le sous-groupe des \'el\'ements de $\Gamma$ qui sont des automorphismes int\'erieurs. On a donc
$$\Gamma_0=\{I(a)/\ a\in \widetilde{K}^*\}\simeq \widetilde{K}^*/Z(M)^*$$
(o\`u, pour un sous-corps $K_0\subset M$ donn\'e, la notation $\widetilde{K_0}$ d\'esigne le commutant de $K_0$ dans $M$).
\vskip 2mm
\noindent
1er point. Puisque $M/L$ est ext\'erieure, on a $\widetilde{L}\subset Z(M)$ et donc $Z(M)=\widetilde{L}$.
\vskip 2mm
\noindent
2\`eme point. L'extension $M/K$ est galoisienne. En effet, le groupe $G_0$ se rel\`eve naturellement dans $\Gamma_0$ et l'on voit imm\'ediatement que $M^{<G_0,G>}=K$.
\vskip 2mm
\noindent
3\`eme point. $<\Gamma_0,G>=\Gamma$. Puisque $M^{<\Gamma_0,G>}=K$, pour montrer que $<\Gamma_0,G>=\Gamma$ il suffit de montrer que $<\Gamma_0,G>$ est un $N$-groupe, c'est-\`a-dire de montrer que l'ensemble
$$A=\{a\in M^*/\ I(a)\in <\Gamma_0,G>\}\cup \{0\}$$
est un corps. Mais, par d\'efinition m\^eme de $\Gamma_0$, on a que $I(a)\in <\Gamma_0,G>\subset \Gamma$ si et seulement si $I(a)\in \Gamma_0$ et donc
$$A=\{a\in M^*/\ I(a)\in \Gamma_0\}\cup \{0\}=\widetilde{K}$$
qui est bien un corps.
\vskip 2mm
\noindent
4\`eme point. $\Gamma_0\triangleleft \Gamma$ et $<\Gamma_0,G>=\Gamma_0\rtimes G$. En effet, si $\sigma\in \Gamma$ et $I(a)\in \Gamma_0$, on a $\sigma\circ I(a)\circ \sigma^{-1}=I(\sigma(a))\in \Gamma_0$. Par ailleurs, si $\sigma=I(a)\in \Gamma_0\cap G$ alors, pour tout $x\in L$, on a $\sigma (x)=axa^{-1}=x$ et donc $a\in \widetilde{L}$. D'apr\`es le point 1, on a $Z(M)=\widetilde{L}$ et l'on en d\'eduit que $\sigma=\hbox{\rm Id}$. L'action de $G$ sur $\Gamma_0$ dans le produit semi-direct $\Gamma_0\rtimes G$ se d\'ecrit en identifiant $\Gamma_0$ au groupe quotient $\widetilde{K}^*/Z(M)^*$ : pour tout $\sigma\in G$ et $\overline{\alpha}\in \widetilde{K}^*/Z(M)^*$, $\overline{\alpha}^\sigma=\sigma(\alpha)$.
\vskip 2mm
\noindent
Conclusion.  Le groupe $\Gamma_0$ est donc un sous-$N$-groupe $N$-invariant de $\Gamma$ et la th\'eorie de Galois assure alors que $M^{\Gamma_0}/K$ est galoisienne de groupe $\Gamma/\Gamma_0=(\Gamma_0\rtimes G)/\Gamma_0\simeq G$, d'apr\`es le point 4. Le dernier point consiste \`a remarquer que
$$x\in M^{\Gamma_0}\Longleftrightarrow \forall a\in \widetilde{K},\ ax=xa\Longleftrightarrow x\in \widetilde{\widetilde{K}}$$
\vskip -2mm
\fin
\vskip 2mm
\noindent
{\bf Corollaire\num\label{2}} (Descente galoisienne) {\it Soit $H$ un corps de dimension finie sur son centre $k$. Si $M/H$ d\'esigne une extension galoisienne de groupe fini $G$, alors l'extension (commutative) $Z(M)/k$ est galoisienne de groupe $G$. Par ailleurs, le corps $M$ est alors isomorphe au produit tensoriel $H\otimes_k Z(M)$.}
\vskip 2mm
\noindent
{\bf Preuve :} Commen\c cons par montrer que l'extension $M/H$ est ext\'erieure en supposant par l'absurde le contraire. Comme $M/H$ est de groupe fini, d'apr\`es \cite[Th\'eor\`eme du \S 2]{Des18}, le centre $C$ de $M$ est un corps fini. Or, $M$ est un corps de dimension (droite ou gauche) finie sur le corps commutatif $k$ et donc, $M$ est aussi de dimension finie sur son centre $C$ (voir par exemple \cite[Lemme 2.1.]{Des01b}). Ainsi, $M$ est lui-m\^eme un corps fini, et donc l'extension $M/H$ est commutative et, par suite, ext\'erieure.
\vskip 2mm
\noindent
\tab On peut donc appliquer le th\'eor\`eme \ref{1}, en posant $K=k$ et $L=H$. L'extension $\widetilde{\widetilde{k}}/k$ est alors galoisienne de groupe $G$. Puisque $Z(H)=k$, on a $H\subset \widetilde{k}$ et donc $\widetilde{\widetilde{k}}\subset \widetilde{H}=Z(M)$ (la derni\`ere \'egalit\'e vient du fait que $M/H$ est ext\'erieure). En reprenant les notations de la preuve du th\'eor\`eme, on a $\widetilde{\widetilde{k}}=M^{\Gamma_0}$. Comme $\Gamma_0$ est compos\'e uniquement d'automorphismes int\'erieurs, on voit que $Z(M)\subset M^{\Gamma_0}$, ce qui prouve finalement que $\widetilde{\widetilde{k}}=Z(M)$.
\vskip 2mm
\noindent
\tab On se retrouve avec les extensions suivantes
$$\xymatrix @!0 @C=3pc @R=3pc{&M&\\
H\ar@{-}[ur]^n_G&&Z(M)\ar@{-}[ul]_{m^2}\\
&k\ar@{-}[ur]^G_n\ar@{-}[ul]^{m^2}&\\}$$
On consid\`ere une $k$-base $e_1,\cdots,e_{m^2}$ de $H\subset M$ que l'on ordonne de sorte $e_1,\cdots ,e_i$ soit une famille $Z(M)$-libre (i.e. lin\'eairement ind\'ependante sur $Z(M)$) maximale pour un certain indice $i$. Si $i\ne m^2$, alors il existe $\lambda_1,\cdots ,\lambda_i\in Z(M)$ tel que $e_{i+1}=\lambda_1e_1+\cdots +\lambda_i e_i$. Puisque tous les $e_j$ sont dans $H$, on a pour tout $\sigma\in G$,
$$e_{i+1}=\sigma(\lambda_1)e_1+\cdots +\sigma(\lambda_i) e_i$$
Par nature m\^eme du centre d'un corps, tout automorphisme de $M$ induit un automorphisme de $Z(M)$ et comme $e_1,\cdots ,e_i$ est $Z(M)$-libre, on en d\'eduit finalement que pour tout $j=1,\cdots ,i$ on a 
$$\forall \sigma\in G,\ \sigma(\lambda_j)=\lambda_j$$
Ceci prouve que $\lambda_j\in k$ pour tout $j=1,\cdots ,i$, ce qui constitue une absurdit\'e, car $e_1,\cdots ,e_{i+1}$ est $k$-libre. Ainsi, la famille $e_1,\cdots,e_{m^2}$ est une $Z(M)$-base de $M$. Par ailleurs, les tenseurs $e_1\otimes 1,\cdots,e_{m^2}\otimes 1$ forment aussi une $Z(M)$-base de $H\otimes_k Z(M)$. L'application $H\otimes_k Z(M)\longrightarrow M$ obtenue par propri\'et\'e universelle du produit tensoriel est donc un isomorphisme de $Z(M)$-alg\`ebres. Ceci prouve bien que $H\otimes_k Z(M)$ est un corps, isomorphe \`a $M$.
\fin
\vskip 2mm
\noindent
\tab Une cons\'equence de ce qui pr\'ec\`ede est alors :
\vskip 2mm
\noindent
{\bf Corollaire\num\label{3}} {\it Pour tout corps $H$ de dimension finie sur son centre $k$, on a $\hbox{\rm PIG}_H\Longrightarrow \hbox{\rm PIG}_k$.}
\vskip 5mm
\noindent
{\bf 1.2.--- Extension des scalaires.}
\vskip 2mm
\noindent
\tab On s'int\'eresse dans cette partie \`a une r\'eciproque du corollaire \ref{2}. On se donne donc un corps $H$ de dimension finie sur son centre $k$ et l'on consid\`ere $L/k$, une extension galoisienne de corps commutatifs de groupe fini $G$. Le corollaire \ref{2} assure que s'il existe une extension galoisienne $M/H$ de groupe fini telle que $L=Z(M)$, alors on a n\'ec\'essairement $M=H\otimes_k L$. On va donc \'etudier ce genre de produits tensoriels, mais pour bien les comprendre, nous allons commencer par rappeler quelques propri\'et\'es sur l'extension des scalaires d'une alg\`ebre simple centrale :
\vskip 2mm
\noindent
\tab On se donne un corps commutatif infini $k$ et une $k$-alg\`ebre simple centrale ${\mathscr A}$. On note $n^2=[{\mathscr A}:k]$ et l'on se donne, une fois pour toute, une $k$-base $(e_1,\cdots ,e_{n^2})$ de ${\mathscr A}$. Si $L/k$ d\'esigne une extension de corps commutatifs (sans hypoth\`ese de finitude ni m\^eme d'alg\'ebricit\'e), l'alg\`ebre obtenue de ${\mathscr A}$ par {\it extension des scalaires \`a $L$} est par d\'efinition la $k$-alg\`ebre tensoris\'ee $\Omega={\mathscr A}\otimes_k L$. On peut plonger $L$ dans $\Omega$ en l'identifiant \`a $k\otimes_k L$. Il est alors clair que $L$ est inclus dans le centre, $Z(\Omega)$, de $\Omega$ de sorte que l'alg\`ebre $\Omega$ peut \^etre consid\'er\'ee comme une $L$-alg\`ebre. Si l'on note $\tilde{e_i}=e_i\otimes 1$ pour tout $i=1,\cdots ,n^2$, on voit que $(\tilde{e_1},\cdots ,\tilde{e_{n^2}})$ est une $L$-base de $\Omega$. On a en particulier $[\Omega:L]=[{\mathscr A}:k]=n^2$. Il est facile de voir qu'en fait, $\Omega$ est une $L$-alg\`ebre simple centrale : le produit tensoriel d'une $k$-alg\`ebre simple centrale par une $k$-alg\`ebre simple est une $k$-alg\`ebre simple (voir par exemple \cite[Th\'eor\`eme II-3]{Bla72}) et le centre d'un produit tensoriel de $k$-alg\`ebres est le produit tensoriel des centres (voir par exemple \cite[Corollaire II-7]{Bla72}). 
\vskip 2mm
\noindent
\tab Rappelons que la norme r\'eduite de ${\mathscr A}/k$ est d\'efinie de la mani\`ere suivante : on commence par se donner un corps neutralisant $D$ de la $k$-alg\`ebre ${\mathscr A}$. Par d\'efinition, il existe un isomorphisme $\varphi : {\mathscr A}\otimes_k D\longrightarrow {\mathscr M}_n(D)$. La norme r\'eduite $\hbox{\rm Nrd}_{{\mathscr A}/k}$ est alors la compos\'ee des applications :
$$\xymatrix{\hbox{\rm Nrd}_{{\mathscr A}/k}:{\mathscr A}\ar[r]^-{a\mapsto a\otimes 1}&{\mathscr A}\otimes_k D\ar[r]^-\varphi&{\mathscr M}_n(D)\ar[r]^-{\hbox{\rm \footnotesize det}}&D}$$
Cette application ne d\'epend, ni du corps neutralisant $D$, ni de l'isomorphisme $\varphi$ et elle est en fait \`a valeurs dans $k$. Une propri\'et\'e importante de la norme r\'eduite est qu'un \'el\'ement $x\in {\mathscr A}$ est inversible si et seulement si $\hbox{\rm Nrd}_{{\mathscr A}/k}(x)\ne 0$ (toutes ces propri\'et\'es sont pr\'esent\'ees dans \cite{Bou12}). L'application ${\cal F}_{{\mathscr A}}(x_1,\cdots ,x_{n^2})\in k[x_1,\cdots ,x_{n^2}]$ d\'efinie pour tout $(x_1,\cdots ,x_{n^2})\in k^{n^2}$, par
$${\cal F}_{\mathscr A}(x_1,\cdots ,x_{n^2})=\hbox{\rm Nrd}_{{\mathscr A}/k}(x_1e_1+\cdots+x_{n^2}e_{n^2})$$
est une $k$-forme de degr\'e $n$ que l'on appellera {\it la forme associ\'ee \`a la norme r\'eduite $\hbox{\rm Nrd}_{{\mathscr A}/k}$ relativement au choix de la $k$-base $(e_1,\cdots ,e_{n^2})$ de ${\mathscr A}$}. 
Par exemple, si ${\mathscr A}={\mathbb H}_k$ d\'esigne l'alg\`ebre des quaternions \`a coefficients dans le corps $k$, alors la norme r\'eduite d'un \'el\'ement $x=a+bu+cv+dw$ (o\`u $u^2=v^2=w^2=uvw=-1$) vaut la traditionnelle norme des quaternions $q(x)=a^2+b^2+c^2+d^2$, de sorte que, pour la choix de la $k$-base $\{1,u,v,w\}$ de ${\mathbb H}_k$, on a
$${\cal F}_{{\mathbb H}_k}(x_1,x_2,x_3,x_4)=x_1^2+x_2^2+x_3^2+x_4^2$$
\vskip 1mm
\noindent
\tab Une des propri\'et\'es importantes que nous allons exploiter dans la suite est donn\'ee par le
\vskip 2mm
\noindent
{\bf Lemme\num\label{5}} {\it Avec les notations pr\'ec\'edentes (on a notamment pos\'e $\Omega={\mathscr A}\otimes_k L$), les formes ${\cal F}_{\mathscr A}$ et ${\cal F}_{\Omega}$ sont \'egales.}
\vskip 2mm
\noindent
{\bf Preuve :} On consid\`ere la cl\^oture alg\'ebrique $\overline{L}$ de $L$. Il existe un $k$-isomorphisme 
$$\varphi :\Omega\otimes_L \overline{L}=\left({\mathscr A}\otimes_k L\right)\otimes_L\overline{L}\longrightarrow {\mathscr A}\otimes_k \overline{L}$$
que l'on peut choisir avec la propri\'et\'e que $\varphi((a\otimes 1)\otimes 1)=a\otimes 1$. Le corps $\overline{L}$ est un corps neutralisant de la $k$-alg\`ebre ${\mathscr A}$ et de la $L$-alg\`ebre $\Omega$, si bien qu'il existe un isomorphisme $\psi: {\mathscr A}\otimes_k \overline{L}\longrightarrow {\mathscr M}_n(\overline{L})$. Le diagramme suivant
$$\xymatrix{\relax {\mathscr A} \ar[r]^-{f}_-{a\mapsto a\otimes 1}\ar[d]^-{\theta}_-{a\mapsto a\otimes 1}&{\mathscr A}\otimes_k \overline{L}\ar[r]^-{\psi}&{\mathscr M}_n(\overline{L})\ar[r]^-{\hbox{\rm \footnotesize det}}&\overline{L}\\
{\mathscr A}\otimes_k L\ar[r]^-{g}_-{x\mapsto x\otimes 1}&({\mathscr A}\otimes_k L)\otimes_L \overline{L}\ar[u]^-{\varphi}\ar[ru]_-{\psi\circ \varphi}}$$
est alors commutatif. Comme $\hbox{\rm Nrd}_{{\mathscr A}/k}=\hbox{\rm \small det}\circ \psi\circ f$ et $\hbox{\rm Nrd}_{\Omega/L}=\hbox{\rm \small det}\circ (\psi\circ \varphi) \circ g$, on en d\'eduit que
$$\hbox{\rm Nrd}_{{\mathscr A}/k}=\hbox{\rm Nrd}_{\Omega/L}\circ \theta$$
Pour tout $(x_1,\cdots ,x_{n^2})\in k^{n^2}$, on a $\theta(x_1e_1+\cdots +x_{n^2}e_{n^2})=x_1\tilde{e_1}+\cdots +x_{n^2}\tilde{e_{n^2}}$ et donc 
$${\cal F}_{\mathscr A}(x_1,\cdots ,x_{n^2})={\cal F}_{\Omega}(x_1,\cdots ,x_{n^2})$$
Puisque le corps $k$ est infini, on en conclut que ${\cal F}_{\mathscr A}={\cal F}_{\Omega}$, par Zariski-densit\'e.
\fin
\vskip 2mm
\noindent
Dans le cas o\`u ${\mathscr A}=H$ d\'esigne un corps, le lemme \ref{5} permet alors de caract\'eriser arithm\'etiquement le cas o\`u $H\otimes_k L$ est un corps :
\vskip 2mm
\noindent
{\bf Proposition\num\label{6}} {\it Soient $H$ un corps de dimension finie sur son centre $k$, ${\cal F}_{H}$ la forme associ\'ee \`a la norme r\'eduite de $H/k$ relativement au choix d'une $k$-base et $L/k$ une extension de corps commutatifs. On a l'\'equivalence 
$$H\otimes_k L\ \hbox{\rm est un corps}\Longleftrightarrow \hbox{\rm la forme ${\cal F}_{H}$ ne poss\`ede que le z\'ero trivial sur $L$}$$}
\vskip 1mm
\noindent
{\bf Preuve :} Comme rappel\'e au d\'ebut de ce paragraphe, l'alg\`ebre $M=H\otimes_k L$ est une $L$-alg\`ebre simple centrale et un \'el\'ement $x\in M$ est inversible si et seulement si on a $\hbox{\rm Nrd}_{M/L}(x)\ne 0$. Puisque d'apr\`es le lemme \ref{5} on a ${\cal F}_{M}={\cal F}_{H}$, l'\'egalit\'e $\hbox{\rm Nrd}_{M/L}(x)= 0$ a lieu pour $x\ne 0$ si et seulement si ${\cal F}_{H}$ poss\`ede un z\'ero non trivial sur $L$.
\fin
\vskip 2mm
\noindent
\tab Revenons maintenant \`a la r\'eciproque du corollaire \ref{2}.
\vskip 2mm
\noindent
{\bf Proposition\num\label{4}} (Extension galoisienne des scalaires) {\it Soit $H$ un corps de dimension finie sur son centre $k$ et $L/k$ une extension galoisienne de corps commutatifs de groupe fini $G$. Si $M=H\otimes_k L$ est un corps, alors l'extension $M/H$ est galoisienne ext\'erieure de groupe $G$.}
\vskip 2mm
\noindent
{\bf Preuve :} On note $\Gamma={\rm Aut}(M/H)$. Puisque $[M:L]=[H:k]=n^2$, on a
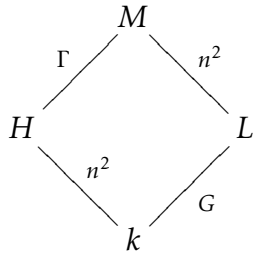
$$\xymatrix @!{&M&\\
H\ar@{-}[ur]^{\Gamma}&&L\ar@{-}[ul]_{n^2}\\
&k\ar@{-}[ur]_G\ar@{-}[ul]_{n^2}&\\}
$$
Pour $\sigma\in {\rm Gal}(L/k)$, on consid\`ere l'application 
$\tilde{\sigma}:H\otimes_k L\longrightarrow H\otimes_k L$ d\'efinie sur les tenseurs $x\otimes t$ par
$$\tilde{\sigma}(x\otimes t)=x\otimes \sigma(t)$$
Puisque $\sigma$ est un $k$-automorphisme de $L$, on voit que $\tilde{\sigma}\in \Gamma$. L'application $\sigma\longmapsto \tilde{\sigma}$ est visiblement un morphisme de groupes et, comme $L$ s'identifie \`a $k\otimes_k L$ dans $H\otimes_k L$, on voit que $\tilde{\sigma}=\hbox{\rm Id}$ seulement si $\sigma=\hbox{\rm Id}$. Ainsi, $G$ se plonge dans $\Gamma$. Il est clair que le sous-corps des \'el\'ements invariants de $M$ par l'action de l'image de $G$ dans $\Gamma$ est \'egal \`a $H\otimes_k k=H$. Il s'ensuit que l'extension $M/H$ est bien galoisienne. 
\vskip 2mm
\noindent
\tab Il reste \`a montrer que $\Gamma$ n'est pas un groupe plus gros que l'image de $G$. Puisque $[M :L]=[H:k]$, par transivit\'e des degr\'es (ici droites ou gauches indistinctement, puisque l'extension $M/H$ est galoisienne) on a que $[M:H]=[L:k]$. Le commutant de $H=H\otimes_k k$ dans $M=H\otimes_k L$ est \'egal \`a $k\otimes_k L=L=Z(M)$ (le commutant d'un produit tensoriel est \'egal au produit tensoriel des commutants, voir par exemple \cite[Proposition II-16]{Bla72}). Il r\'esulte de ce dernier fait que l'extension $M/H$ est ext\'erieure et, en vertu des propri\'et\'es g\'en\'erales de la th\'eorie de Galois des corps gauches rappel\'ees dans le \S 1.1., on a donc
$$|\Gamma|=[M:H]=[L:k]=|G|$$
Puisque $G$ se plonge dans $\Gamma$, on en d\'eduit finalement que $\Gamma=G$. 
\fin
\vskip 2mm
\noindent
\tab On voit donc que r\'ealiser $G$ sur $H$ revient juste \`a r\'ealiser $G$ par une extension de corps commutatifs $L/k$ de sorte que $H\otimes_k L$ soit un corps. Par exemple, si l'on reprend l'alg\`ebre de quaternions ${\mathbb H}_k$, alors la proposition \ref{6} montre que ${\mathbb H}_k\otimes_kL$ est un corps, si et seulement si $L$ est de niveau\footnote{Rappelons que le niveau, $\nu(L)$, d'un corps commutatif $L$ est, $+\infty$ si $-1$ n'est pas une somme de carr\'ees dans $L$ (la th\'eorie d'Artin-Schreier montre que cette propri\'et\'e \'equivaut au fait d'\^etre ordonnable pour $L$, voir \cite{Rib72}), et dans le cas contraire, le niveau est le plus petit entier $n$ tel que $-1$ soit la somme de $n$ carr\'es dans $L$.} $\nu(L)\geq 4$. En effet, dire que la forme ${\cal F}_{{\mathbb H}_k}(x_1,x_2,x_3,x_4)=x_1^2+x_2^2+x_3^2+x_4^2$ poss\`ede un z\'ero non trivial sur $L$ \'equivaut \`a dire que $-1$ est somme de trois carr\'es dans $L$. On voit, en particulier, que ${\mathbb H}_k\otimes_kL$ est un corps d\`es que $L$ est ordonnable. Puisque $k={\mathbb R}$ est un corps r\'eel clos on en d\'eduit que le corps des quaternion d'Hamilton ${\mathbb H}={\mathbb H}_{\hbox{\scriptsize $\mathbb R$}}$ ne poss\`ede aucun extension galoisienne finie non triviale (en fait ${\mathbb H}$ ne poss\`ede aucune extension finie non triviale, voir \cite{Des01b}).
\vskip 2mm
\noindent
\tab En application du corollaire \ref{2} et des propositions \ref{6} et \ref{4}, on d\'eduit le
\vskip 2mm
\noindent
{\bf Th\'eor\`eme\num\label{7}} {\it Soient $H$ un corps de dimension finie sur son centre $k$, ${\cal F}_{H}$ la forme associ\'ee \`a la norme r\'eduite de $H/k$ relativement au choix d'une $k$-base de $H$ et $G$ un groupe fini. Pour qu'il existe une extension galoisienne $M/H$ de groupe $G$ il faut et il suffit qu'il existe une extension galoisienne $L/k$ de groupe $G$ telle que ${\cal F}_{H}$ ne poss\`ede que le z\'ero trivial sur $L$. Dans ces conditions, on a alors $Z(M)=L$ et $M=H\otimes_k L$.}
\vskip 2mm
\noindent
d'o\`u d\'ecoule le th\'eor\`eme A \'enonc\'e dans l'introduction.
\vskip 1cm
\noindent
{\bf \large 2.--- Applications aux corps des fractions tordus.}
\vskip 2mm
\noindent
\tab On souhaite maintenant appliquer ce qui pr\'ec\`ede au cas des corps de fractions rationnelles tordus \`a ind\'etermin\'ee centrale et \`a coefficients dans $H$ lorsque $H$ est un corps de dimension finie sur son centre $k$. Ce cas rentre parfaitement dans le cadre d'\'etude du paragraphe 1 puisque, comme nous allons le voir plus bas, $H(t)$ est une $k(t)$-alg\`ebre simple centrale. Etablissons pr\'ealablement un lemme utile pour la suite de ce texte :
\vskip 2mm
\noindent
{\bf Lemme\num\label{lem}} {\it Soient $k$ un corps commutatif et ${\cal F}$ une forme \`a coefficients dans $k$. Si ${\cal F}$ ne poss\`ede que le z\'ero trivial sur $k$, alors il en est de m\^eme sur le corps $k((t))$ des s\'eries de Laurent \`a coefficients dans $k$.}
\vskip 2mm
\noindent
{\bf Preuve :} Supposons que ${\cal F}\in k[x_1,\cdots ,x_d]$ poss\`ede un z\'ero non trivial $(r_1(t),\dots ,r_{d}(t))\in k((t))^d$. Quitte \`a factoriser par la puissance de $t$ correspondant \`a la plus petite des valuations des s\'eries $r_i(t)$, on peut supposer que $(r_1(t),\dots ,r_{d}(t))\in k[[t]]^{d}$ et qu'au moins un des $r_i(t)$ est de valuation nulle. En posant $t=0$, on voit que $(r_1(0),\dots ,r_{d}(0))\in k^{d}$ est alors un z\'ero non trivial de ${\cal F}$.
\fin
\vskip 2mm
\noindent
{\bf Proposition\num\label{8}} {\it L'alg\`ebre $H\otimes_k k(t)$ est un corps, isomorphe au corps $H(t)$ des fractions rationnelles tordu \`a ind\'etermin\'ee centrale.}
\vskip 2mm
\noindent
{\bf Preuve :} On reprend les notations du \S 1.2. L'alg\`ebre $H\otimes_k k(t)$ est une $k(t)$-alg\`ebre simple centrale et le lemme \ref{5} assure que l'on a ${\cal F}_{H\otimes_k k(t)}={\cal F}_{H}$. Puisque $H$ est un corps, la forme ${\cal F}_{H}$ ne poss\`ede que le z\'ero trivial sur $k$ et le lemme \ref{lem} prouve qu'il en est alors de m\^eme sur $k(t)$. Il d\'ecoule alors de la proposition \ref{6} que l'alg\`ebre $H\otimes_k k(t)$ est bien un corps. 
\vskip 2mm
\noindent
\tab Par propri\'et\'e universelle du produit tensoriel, l'application 
$$(a,r(t))\in H\times k(t)\longmapsto a.r(t)\in H(t)$$
d\'efinit un morphisme de $k$-espaces vectoriels $\psi:H\otimes_k k(t)\longrightarrow H(t)$. Par d\'efinition, dans $H(t)$, on a $a.r(t)=r(t).a$ et l'on voit donc que $\psi$ est un morphisme de $k$-alg\`ebres. Ce dernier est alors injectif puisque $H\otimes_k k(t)$ est un corps. Il est clair que $\psi$ identifie les anneaux $H\otimes_k k[t]$ et $H[t]$. Comme d\'ej\`a rappel\'e dans l'introduction, les \'el\'ements de $H(t)$ s'\'ecrivent sous la forme $p(t)q(t)^{-1}$ avec $p(t),q(t)\in H[t]$. Puisque $H\otimes_k k(t)$ est un corps, on en d\'eduit finalement que $\psi$ est aussi surjectif.
\fin
\vskip 2mm
\noindent
\tab Avec les notations pr\'ec\'edentes, on a donc $\mathcal{F}_{H(t)}=\mathcal{F}_H$ et, par application du th\'eor\`eme A, on en d\'eduit que ${\rm{PIG}}{\mathcal{F}_H}{\rm{C}}_{k(t)} \Longleftrightarrow \hbox{\rm PIG}_{H(t)}$. Le lemme \ref{lem} montre que, pour toute $k$-forme $\mathcal{F}$ ne poss\'edant que le z\'ero trivial sur $k$, on a l'implication
$$(1)\ \ \ \hbox{PIGFR}_{k} \Longrightarrow {\rm{PIG}}{\mathcal{F}}{\rm{C}}_{k(t)}$$
o\`u $\hbox{PIGFR}_{k}$ d\'esigne la variante du PIG$_{k(t)}$ suivante :
\vskip 2mm
\noindent
{\bf $\hbox{\bf PIGFR}_k$ (Probl\`eme Inverse de Galois Fortement R\'egulier)} : {\it tout groupe fini $G$ est-il groupe de Galois d'une extension galoisienne $L/k(t)$ telle que $L$ se plonge dans le corps $k((t))$ ?}
\vskip 2mm
\noindent
\tab D'un point de vue g\'eom\'etrique, le $\hbox{\rm PIGFR}_k$ \'equivaut \`a savoir si, pour tout groupe fini $G$, il existe un rev\^etement galoisien $X\longrightarrow {\mathbb P}^1$ de groupe $G$, d\'efini sur $k$ et tel que la courbe $X$ soit g\'eom\'etriquement irr\'eductible et poss\`ede un point $k$-rationnel non ramifi\'e. Lorsque $k$ est ample, pour $X\longrightarrow {\mathbb P}^1$ donn\'e, l'existence d'un point $k$-rationnel non ramifi\'e est en fait \'equivalente \`a l'existence d'au moins un point $k$-rationnel. Si l'on retire \`a cette formulation g\'eom\'etrique du $\hbox{\rm PIGFR}_k$ l'hypoth\`ese d'existence d'un point $k$-rationnel, on retombe sur le classique
\vskip 2mm
\noindent
{\bf $\hbox{\bf PIGR}_k$ (Probl\`eme Inverse de Galois R\'egulier)} : {\it tout groupe fini est-il groupe de Galois d'une extension galoisienne $L/k(t)$ telle que $L \cap \overline{k}=k$ ?}
\vskip2mm
\noindent
et l'on a donc l'implication 
$$(2)\ \ \ {\rm{PIGFR}}_k \Longrightarrow {\rm{PIGR}}_k$$
ce qui justifie la terminologie "fortement r\'egulier" pour le $\hbox{\rm PIGFR}$. Il est conjectur\'e que le PIGR admet une r\'eponse positive pour tout corps mais personne ne s'est encore aventur\'e \`a conjecturer qu'il en \'etait de m\^eme pour le \hbox{\rm PIGFR}. La litt\'erature ne contient, \`a l'heure actuelle, aucun exemple de corps commutatif ayant une r\'eponse n\'egative \`a ce probl\`eme.
\vskip 2mm
\noindent
{\bf Proposition\num\label{9}} {\it Pour toute extension de corps commutatifs $L/k$, on a l'implication
$${\rm{PIGFR}}_{k} \Longrightarrow {\rm{PIGFR}}_{L}$$
En particulier, si le corps $k$ contient un corps ample alors le ${\rm{PIGFR}}_{k}$ admet une r\'eponse positive.}
\vskip 1mm
\noindent
{\bf Preuve :} Fixons un groupe fini $G$ et un \'el\'ement transcendant $t$ sur $L$. Par hypoth\`ese, il existe une extension galoisienne $M/k(t)$ de groupe $G$ telle que $M$ se plonge dans $k((t))\subset L((t))$. Cette derni\`ere condition entra\^ine $M \cap \overline{k}=k$. L'extension galoisienne $ML/L(t)$ est alors galoisienne de groupe $G$, m\^eme si $L/k$ n'est pas n\'ecessairement alg\'ebrique (voir, par exemple, \cite[Lemma 16.2.1]{FJ08} ou \cite[Proposition 2.3.2]{Deb09}), et l'on a $ML \subseteq L((t))$.
\vskip 2mm
\noindent
\tab Le fait que le ${\rm{PIGFR}}_{k}$ admette alors une r\'eponse positive lorsque $k$ est ample fut d\'emontr\'e par Pop dans \cite{Pop96} (voir aussi \cite[Theorem C]{HJ98}). Ce r\'esultat prouve la fin de la proposition.
\fin
\vskip 2mm
\noindent
\tab Le th\'eor\`eme B de l'introduction d\'ecoule finalement de tous ces r\'esultats : si le corps $k$ contient un corps ample, alors la proposition \ref{9} assure que le ${\rm{PIGFR}}_{k}$ admet une r\'eponse positive et l'implication (1) montre qu'il en est de m\^eme du ${\rm{PIG}}{\mathcal{F}_H}{\rm{C}}_{k(t)}$. Le th\'eor\`eme A permet alors de conclure.
\vskip 2mm
\noindent
\tab Pour compl\'eter les liens implicatifs que nous avons \'etablis, notons le r\'esultat suivant (qui est trivial dans le cas des corps commutatifs) :
\vskip 2mm
\noindent
{\bf Proposition\num\label{10}} {\it Si $H$ d\'esigne un corps gauche de dimension finie sur son centre $k$ alors $
{\rm{PIG}}_H \Longrightarrow {\rm{PIG}}_{H(t)}$.}
\vskip 2mm
\noindent
{\bf Preuve :} Consid\'erons une $k$-forme $\mathcal{F}$ ne poss\'edant que le z\'ero trivial sur $k$ et un groupe fini $G$. S'il existe une extension galoisienne $L/k$ de groupe $G$ telle que $\mathcal{F}$ ne poss\`ede que le z\'ero trivial sur $L$ alors l'extension $L(t)/k(t)$ est galoisienne de groupe $G$ et le lemme \ref{lem} montre que $\mathcal{F}$ ne poss\`ede que le z\'ero trivial sur $L(t)$. On a donc l'implication 
$$(3)\ \ \ {\rm{PIG}}{\mathcal{F}}{\rm{C}}_{k} \Longrightarrow {\rm{PIG}}{\mathcal{F}}{\rm{C}}_{k(t)}$$
La proposition d\'ecoule alors du th\'eor\`eme A en exploitant l'implication ci-dessus pour le choix de $\mathcal{F}=\mathcal{F}_H=\mathcal{F}_{H(t)}$.
\fin
\vskip 2mm
\noindent
\tab Le diagramme d'implications suivant r\'esume les r\'esultats \'etablis dans ce texte :
\vskip 6mm
\noindent
$$\footnotesize \xymatrix@!0 @R=6.5pc @C=8.5pc{\relax
&&&\txt{\footnotesize $\hbox{\rm PIG${\cal F}$C}_k$\\ \scriptsize cas a (resp. cas b)}\ar@{=>}[dl]_-{(3)}\ar@{=>}[d]^*[@]{\hbox to 2pt{\hss\txt{\scriptsize (Th.A)}\hss}}\ar@{=>}[rd]\ar@/_5mm/@{<=}[d]_*[@]{\hbox to 9pt{\hss\txt{\scriptsize cas b}\hss}}&\\
&&\txt{\footnotesize $\hbox{\rm PIG${\cal F}$C}_{k(t)}$\\ \scriptsize cas a (resp. cas b)}\ar@/_7mm/@{<=}[dr]_*[@]{\hbox to 3pt{\hss\txt{\scriptsize cas b}\hss}}\ar@{=>}[dr]^*[@]{\hbox to 2pt{\hss\txt{\scriptsize (Th.A)}\hss}}_*[@]{\hbox to 3pt{\hss\txt{\scriptsize cas a et b}\hss}}&\ \ \txt{\footnotesize $\hbox{\rm PIG}_{H}$\\ \scriptsize cas a (resp. cas b)}\ar@{=>}[d]^*[@]{\hbox to 2pt{\hss\txt{\scriptsize (Prop.\ref{10})}\hss}}\ar@{=>}[r]^-{\hbox{\scriptsize (Coro.\ref{3})}}_-{\hbox{\scriptsize cas b}}&\hbox{\rm PIG}_k\ar@/_5mm/@{=>}[ul]_{\txt{\scriptsize cas a}}\\
\txt{\footnotesize $k$ contient un\\ \footnotesize corps ample}\ \ \ar@{=>}[r]^-{(\hbox{\scriptsize Prop.} \ref{9})} &\hbox{\rm PIGFR}_k\ar@{=>}[rd]\ar@{=>}[ru]^-{(1)}\ar@{=>}[rd]^-{(2)}&&\ \ \txt{\footnotesize $\hbox{\rm PIG}_{H(t)}$\\ \scriptsize cas a (resp. cas b)}&\\
&&\hbox{\rm PIGR}_k\ar@{=>}[ru]_*[@]{\hbox to 0pt{\hss\txt{\scriptsize cas a}\hss}}&&\\}$$
\vskip 4mm
\noindent
{\footnotesize
Cas a (resp. Cas b) : $H$ est un corps commutatif (resp. gauche) de dimension finie sur $k$ (resp. sur son centre $k$) et ${\cal F}$ d\'esigne la forme associ\'ee \`a la norme (resp. la norme r\'eduite) de $H/k$ relativement au choix d'une $k$-base de $H$.}
\vskip 2mm
\noindent
\tab Nous donnons enfin quelques exemples de groupes finis se r\'ealisant comme groupes de Galois sur $H(t)$ sans supposer que le centre de $H$ ne contienne de corps ample :
\vskip 2mm
\noindent
{\bf{Proposition \num\label{11}}} {\it Soient $G$ un groupe fini et $H$ un corps de dimension finie sur son centre $k$. Dans chacun des cas suivants, $G$ est le groupe de Galois d'une extension galoisienne de $H(t)$ : 
\vskip 2mm
\noindent
{\rm{1/}} $G$ est ab\'elien,
\vskip 2mm
\noindent
{\rm{2/}} $G=S_n$ ($n \geq 3$) et $k$ est infini,
\vskip 2mm
\noindent
{\rm{3/}} $G=A_n$ ($n \geq 4$) et $k$ est de caract\'eristique nulle,
\vskip 2mm
\noindent
{\rm{4/}} $G$ est r\'esoluble et $k$ est de caract\'eristique strictement positive.}
\vskip 2mm
\noindent 
{\bf Preuve :} Il suffit, dans chaque cas, de trouver un sous-corps $k_0$ de $k$ et une extension galoisienne $L/k_0(t)$ de groupe $G$ telle que $L$ se plonge dans le corps $k_0((t))$. L'extension $Lk/k(t)$ est alors une galoisienne de groupe $G$ qui se plonge dans $k((t))$. Cette hypoth\`ese de plongement permet alors, gr\^ace au lemme \ref{lem}, d'appliquer le th\'eor\`eme \ref{7} au corps $H(t)$.
\vskip 2mm
\noindent
1/ Si $G$ est ab\'elien, cette propri\'et\'e est vraie sur le corps $k$, en tant que cons\'equence bien connue du {\it{twisting lemma}} (voir \cite{Deb99a} et \cite[Proposition 3.2.4]{Deb09}) et de l'existence d'une extension galoisienne $L/k(t)$ de groupe $G$ telle que l'id\'eal $\langle t \rangle$ de $\overline{k}[t]$ ne se ramifie pas dans l'extension $L\overline{k}/\overline{k}(t)$ (voir \cite{Deb99e}).
\vskip 2mm
\noindent
2/ Puisque $k$ est infini, il poss\`ede un sous-corps $k_0$ qui est soit ample, soit hilbertien. En effet, si $k$ est de caract\'eristique nulle, on peut prendre $k_0=\mathbb{Q}$. Si $k$ est de caract\'eristique $p>0$ il est alors, soit extension alg\'ebrique infinie du corps $\mathbb{F}_p$ et donc $k_0=k$ est un corps PAC (voir \cite[Corollary 11.2.4]{FJ08}), soit une extension du corps de fractions rationnelles $k_0={\mathbb F}_p(x)$ qui est hilbertien. Le cas ample ayant d\'ej\`a \'et\'e trait\'e de mani\`ere plus g\'en\'erale, on  suppose d\'esormais que $k_0$ est hilbertien. On se donne alors un polyn\^ome unitaire s\'eparable $P_0(x) \in k_0[x]$ de degr\'e $n$ et dont toutes les racines sont dans $k_0$, et un polyn\^ome unitaire $P_1(x) \in k_0[x]$ de degr\'e $n$ et de groupe de Galois $S_n$ sur $k_0$ (un tel polyn\^ome existe par \cite[Corollary 16.2.7]{FJ08}). Par interpolation polynomiale, il existe un polyn\^ome unitaire $P(t,x) \in k_0[t][x]$ de degr\'e $n$ tel que $P(0,x)=P_0(x)$ et $P(1,x) = P_1(x)$. Notons $L$ le corps de d\'ecomposition sur $k_0(t)$ de $P(t,x)$. Puisque $P(0,x)$ est s\'eparable et a toutes ses racines dans $k_0$, le corps $L$ se plonge dans $k_0((t))$. Enfin, puisque $P(1,x)$ poss\`ede $S_n$ comme groupe de Galois sur $k_0$, le groupe de Galois de $L/k_0(t)$ vaut \'egalement $S_n$.
\vskip 2mm
\noindent
3/ Ici $k=k_0$ convient. En effet, fixons un polyn\^ome unitaire s\'eparable $P(x)\in k[x]$ de degr\'e $n$ et dont toutes les racines sont dans $k$. Par \cite[Theorem 3]{KM01}, il existe un polyn\^ome unitaire $P(t,x)\in k[t][x]$ tel que, si $L$ est le corps de d\'ecomposition sur $k(t)$ de $P(t,x)$, alors ${\rm{Gal}}(L/k(t))=A_n$ et le corps de d\'ecomposition sur $k$ de $P(0,x)$ est \'egal \`a celui de $P(x)$, c'est-\`a-dire \`a $k$. Comme $P(x)$ est s\'eparable, cela entra\^ine que $L$ se plonge dans $k((t))$.
\vskip 2mm
\noindent
4/ La propri\'et\'e est vraie pour $k_0$ \'egal au sous-corps premier de $k$, en tant que cons\'equence de \cite[p. 597, Exercise (a)]{NSW08}.
\fin
\vskip 2mm
\noindent
\tab Pour conclure cet article, soulignons que si PIG est tout \`a fait l\'egitime pour les corps gauches, le probl\`eme r\'egulier et toutes les autres variantes que nous avons introduites dans ce texte n'ont pas de sens {\it a priori}. Il faut avoir \`a l'esprit que les notions d'\'el\'ements alg\'ebriques et de cl\^oture alg\'ebrique n'ont rien de naturel pour un corps gauche. Ce point rend beaucoup plus d\'elicat l'approche au cas non commutatif de la th\'eorie de Galois. A titre d'exemple, il est ainsi montr\'e dans \cite{Des01b} que les extensions finies d'un corps alg\'ebriquement clos $\overline{k}$ de caract\'eristique $0$ sont toutes de degr\'e $2$ et qu'il en existe une infinit\'e non isomorphes deux \`a deux. Dans cette situation, si l'on plonge deux telles extensions dans un corps, alors le compositum est n\'ecessairement de degr\'e infini sur $\overline{k}$.

\vskip 1cm
\bibliography{Biblio2}

\begin{thebibliography}{NSW08}

\bibitem[Bla72]{Bla72}
Andr\'e Blanchard.
\newblock {\em Les corps non commutatifs. ({F}rench)}.
\newblock Collection {S}up : {L}e {M}ath\'ematicien, {N}o. 9. Presses
  {U}niversitaires de {F}rance, {V}end\^ome, 1972.
\newblock 135 pp.

\bibitem[Bou12]{Bou12}
Nicolas Bourbaki.
\newblock {\em \'{E}l\'ements de math\'ematique. {A}lg\`ebre. {C}hapitre 8.
  {M}odules et anneaux semi-simples. ({F}rench)}.
\newblock Springer, Berlin, 2012.
\newblock x+489 pp. Second revised version of the 1958 edition.

\bibitem[BSF13]{BSF13}
Lior Bary-Soroker and Arno Fehm.
\newblock Open problems in the theory of ample fields.
\newblock In {\em Geometric and differential {G}alois theories}, volume~27 of
  {\em S\'emin. Congr.}, pages 1--11. Soc. Math. France, Paris, 2013.

\bibitem[Coh95]{Coh95}
Paul~Moritz Cohn.
\newblock {\em Skew fields. {T}heory of general division rings}.
\newblock Encyclopedia of {M}athematics and its {A}pplications, 57. Cambridge
  {U}niversity {P}ress, Cambridge, 1995.
\newblock xvi + 500 pp.

\bibitem[DD97]{DD97b}
Pierre D{\`e}bes and Bruno Deschamps.
\newblock The regular inverse {G}alois problem over large fields.
\newblock In {\em Geometric {G}alois actions, 2}, volume 243 of {\em London
  Math. Soc. Lecture Note Ser.}, pages 119--138. Cambridge Univ. Press,
  Cambridge, 1997.

\bibitem[D{\`{e}}b99a]{Deb99a}
Pierre D{\`{e}}bes.
\newblock Galois covers with prescribed fibers: the {B}eckmann-{B}lack problem.
\newblock {\em Ann. Scuola Norm. Sup. Pisa Cl. Sci. (4)}, 28(2):273--286, 1999.

\bibitem[D{\`e}b99b]{Deb99e}
Pierre D{\`e}bes.
\newblock Regular realization of abelian groups with controlled ramification.
\newblock In {\em Applications of curves over finite fields (Seattle, WA,
  1997)}, volume 245 of {\em Contemp. Math.}, pages 109--115. Amer. Math. Soc.,
  Providence, RI, 1999.

\bibitem[D{\`e}b09]{Deb09}
Pierre D{\`e}bes.
\newblock {\em Arithm\'etique des rev\^etements de la droite}.
\newblock {L}ecture notes, 2009.
\newblock {A}t http://math.univ-lille1.fr/\~{}pde/ens.html.

\bibitem[Des01]{Des01b}
Bruno Deschamps.
\newblock A propos d'un th\'eor\`eme de {F}robenius. ({F}rench).
\newblock {\em Ann. {M}ath. {B}laise {P}ascal}, 8(2):61--66, 2001.

\bibitem[Des18]{Des18}
Bruno Deschamps.
\newblock Des extensions plus petites que leurs groupes de {G}alois.
  ({F}rench).
\newblock {\em Comm. Algebra}, 46(10):4555--4560, 2018.

\bibitem[FHV93]{FHV93}
Michael~D. Fried, Dan Haran, and Helmut V\"{o}lklein.
\newblock Absolute {G}alois group of the totally real numbers.
\newblock {\em C. R. Acad. Sci. Paris S\'er. I Math.}, 317(11):995--999, 1993.

\bibitem[FJ08]{FJ08}
Michael~D. Fried and Moshe Jarden.
\newblock {\em Field arithmetic}.
\newblock Ergebnisse der Mathematik und ihrer Grenzgebiete. 3. Folge. A
  {S}eries of Modern Surveys in Mathematics [Results in Mathematics and Related
  Areas. 3rd Series. A Series of Modern Surveys in Mathematics], 11.
  Springer-Verlag, Berlin, third edition, 2008.
\newblock Revised by Jarden. xxiv + 792 pp.

\bibitem[GJ02]{GJ02}
Wulf-Dieter Geyer and Moshe Jarden.
\newblock {P}{{\it{{S}}}}{C} {G}alois extensions of {H}ilbertian fields.
\newblock {\em Math. Nachr.}, 236:119--160, 2002.

\bibitem[HJ98]{HJ98}
Dan Haran and Moshe Jarden.
\newblock Regular split embeddings problems over complete valued fields.
\newblock {\em Forum Math.}, 10(3):329--351, 1998.

\bibitem[KM01]{KM01}
J\"urgen Kl\"uners and Gunter Malle.
\newblock A database for field extensions of the rationals.
\newblock {\em LMS J. Comput. Math.}, 4:182--196, 2001.

\bibitem[NSW08]{NSW08}
J\"urgen Neukirch, Alexander Schmidt, and Kay Wingberg.
\newblock {\em Cohomology of number fields}, volume 323 of {\em Grundlehren der
  mathematischen Wissenschaften [{F}undamental {P}rinciples of {M}athematical
  {S}ciences]}.
\newblock Springer-Verlag, Berlin, second edition, 2008.
\newblock xvi+825 pp.

\bibitem[Ore33]{Ore33}
Oystein Ore.
\newblock Theory of non-commutative polynomials.
\newblock {\em Ann. of Math. (2)}, 34(3):480--508, 1933.

\bibitem[Pop96]{Pop96}
Florian Pop.
\newblock Embedding problems over large fields.
\newblock {\em Ann. of Math. (2)}, 144(1):1--34, 1996.

\bibitem[Pop14]{Pop14}
Florian Pop.
\newblock Little survey on large fields - old $\&$ new.
\newblock In {\em Valuation theory in interaction}, EMS Ser. Congr. Rep., pages
  432--463. Eur. Math. Soc., Z\"urich, 2014.

\bibitem[Rib72]{Rib72}
Paulo Ribenboim.
\newblock {\em L'arithm\'etique des corps. ({F}rench)}.
\newblock Hermann, Paris, 1972.
\newblock 245 pp.

\end{thebibliography}
\bibliographystyle{alpha}
\vskip 8mm
\noindent
{\bf Bruno Deschamps}\\
{\sc Laboratoire de Math\'ematiques Nicolas Oresme, CNRS UMR 6139}\\
Universit\'e de Caen - Normandie\\
BP 5186, 14032 Caen Cedex - France
\vskip -1.5mm
\noindent
------------------------------
\vskip -1.5mm
\noindent
{\sc D\'epartement de Math\'ematiques --- Le Mans Universit\'e}\\
Avenue Olivier Messiaen, 72085 Le Mans cedex 9 - France\\
E-mail : Bruno.Deschamps@univ-lemans.fr
\vskip 4mm
\noindent
{\bf Fran\c cois Legrand}\\
{\sc Institut f\"ur Algebra, Fachrichtung Mathematik}\\
TU Dresden, 01062 Dresden, Germany\\
E-mail : francois.legrand@tu-dresden.de

\end{document}